\documentclass[a4paper,12pt,reqno]{amsart}
\usepackage{amsmath}
\usepackage{amssymb}
\usepackage{amsthm}

\usepackage{amscd}
\usepackage{amsfonts}
\usepackage{setspace}
\usepackage{version}

\title[Moments of the zeta function]{Sharp conditional bounds for moments of the Riemann zeta function}
\author{Adam J Harper}
\address{Centre de recherches math\'{e}matiques, Universit\'{e} de Montr\'{e}al, Pavillon Andr\'{e}-Aisenstadt, 2920 Chemin de la tour, bur. 5357, Montr\'{e}al QC H3T 1J4, Canada}
\email{harperad@crm.umontreal.ca}
\date{20th May 2013}
\thanks{The author is supported by a postdoctoral fellowship from the Centre de recherches math\'{e}matiques in Montr\'{e}al.}

\numberwithin{equation}{section}
\theoremstyle{plain}

\onehalfspace
\newcommand{\N}{\mathbb{N}}

\newcommand{\E}{\mathbb{E}}

\newcommand{\C}{\mathbb{C}}

\newtheorem{thm1}{Theorem}
\newtheorem{prop1}{Proposition}
\newtheorem{prop2}[prop1]{Proposition}
\newtheorem{lem1}{Lemma}
\newtheorem{lem2}[lem1]{Lemma}
\newtheorem{lem3}[lem1]{Lemma}
\newtheorem{thm2}[thm1]{Theorem}
\newtheorem{prop3}[prop1]{Proposition}
\newtheorem{prop4}[prop1]{Proposition}

\begin{document}

\maketitle

\begin{abstract}
We prove, assuming the Riemann Hypothesis, that $\int_{T}^{2T} |\zeta(1/2+it)|^{2k} dt \ll_{k} T \log^{k^{2}}T$ for any fixed $k \geq 0$ and all large $T$. This is sharp up to the value of the implicit constant.

Our proof builds on well known work of Soundararajan, who showed, assuming the Riemann Hypothesis, that $\int_{T}^{2T} |\zeta(1/2+it)|^{2k} dt \ll_{k, \epsilon} T \log^{k^{2}+\epsilon}T$ for any fixed $k \geq 0$ and $\epsilon > 0$. Whereas Soundararajan bounded $\log|\zeta(1/2+it)|$ by a single Dirichlet polynomial, and investigated how often it attains large values, we bound $\log|\zeta(1/2+it)|$ by a sum of many Dirichlet polynomials and investigate the joint behaviour of all of them. We also work directly with moments throughout, rather than passing through estimates for large values.
\end{abstract}

\section{Introduction}
In this paper we shall prove the following result.

\begin{thm1}
Assume the Riemann Hypothesis is true, and let $k \geq 0$ be fixed. Then for all large $T$ we have
$$ \int_{T}^{2T} |\zeta(1/2+it)|^{2k} dt \ll_{k} T \log^{k^{2}}T , $$
where the implicit constant depends on $k$ only.
\end{thm1}

Theorem 1 is sharp up to the value of the implicit constant, since Radziwi\l{\l} and Soundararajan~\cite{radsound} have proved a matching lower bound unconditionally when $k \geq 1$, and earlier Ramachandra (see \cite{ramachandra, ramachandra2}) and Heath-Brown~\cite{hbfrac} proved a matching lower bound for all $k \geq 0$, but assuming the Riemann Hypothesis. See e.g. the introduction to Soundararajan's paper~\cite{soundmoments} for further references.

The theorem is known unconditionally when $k=0,1,2$, due to classical work of Hardy--Littlewood and Ingham, and has been known for a long time, assuming the Riemann Hypothesis, for all $k \leq 2$, due to work of Ramachandra and Heath-Brown~\cite{hbfrac}.

Recently Radziwi\l\l~\cite{rad1} proved the theorem on the range $2 < k < 2 + 2/11$, conditionally on the Riemann Hypothesis, by a nice argument using an estimate for the integral of $|\zeta(1/2+it)|^{4}$ multiplied by the square of a short Dirichlet polynomial. Our proof of Theorem 1 is not like this, but instead extends an argument of Soundararajan~\cite{soundmoments}, who showed, assuming the Riemann Hypothesis, that
$$ \int_{T}^{2T} |\zeta(1/2+it)|^{2k} dt \ll_{k, \epsilon} T \log^{k^{2}+\epsilon}T  $$
for any fixed $k \geq 0$ and $\epsilon > 0$. This estimate was itself a considerable improvement over the conditional bound $\int_{T}^{2T} |\zeta(1/2+it)|^{2k} dt \ll T e^{O(k\log T /\log\log T)}$, which follows because the Riemann Hypothesis implies the pointwise bound $|\zeta(1/2+it)| \leq e^{O(\log T / \log\log T)}$, for $T \leq t \leq 2T$.

\vspace{12pt}
To explain the proof of Theorem 1, we first give a brief heuristic discussion about the behaviour of $\log|\zeta(1/2+it)|$. Note that $\int_{T}^{2T} |\zeta(1/2+it)|^{2k} dt = \int_{T}^{2T} e^{2k\log|\zeta(1/2+it)|} dt$, so if we can understand how $\log|\zeta(1/2+it)|$ behaves as $t$ varies then we can understand the moments of the zeta function.

When $\Re(s) > 1$ the zeta function is given by an absolutely convergent Euler product, and so $\log|\zeta(s)|$ is given by an absolutely convergent sum over primes. On the critical line $\Re(s)=1/2$ this is no longer true, since the zeros of the zeta function have an influence: but, assuming the Riemann Hypothesis, the influence of the zeros can be quite well understood, and $\zeta(1/2+it)$ still looks ``quite a lot'' like an Euler product. For example, Gonek, Hughes and Keating~\cite{gonekhugheskeating} showed, roughly speaking, that if the Riemann Hypothesis is true then
$$ \zeta(1/2+it) \approx \prod_{p \leq T} \left(1- \frac{1}{p^{1/2+it}}\right)^{-1} \cdot \prod_{|\gamma - t| < 1/\log T} (ci(t-\gamma)\log T), $$
where $p$ denotes primes, $c > 0$ is a constant, and the second product is over ordinates $\gamma$ of zeros of the zeta function. Thus we have
$$ \log|\zeta(1/2+it)| = \Re \log\zeta(1/2+it) \approx \Re \sum_{p \leq T} \frac{1}{p^{1/2+it}} + \text{contribution from zeros}. $$
Here we note two things about the ``contribution from zeros'':
\begin{enumerate}
\item For most values $T \leq t \leq 2T$ this contribution will have size $O(1)$, since the averaging spacing between zeros at height $t$ is $\Theta(1/\log T)$. Such a contribution can be absorbed into the implicit constant in moment bounds, which provides hope that one can sharply understand the moments of zeta by understanding Dirichlet polynomials like $\sum_{p \leq T} \frac{1}{p^{1/2+it}}$.

\item For {\em all} values $T \leq t \leq 2T$ this contribution is negative (roughly speaking), so to obtain upper bounds for moments one {\em only} needs to consider the Dirichlet polynomial contribution. This kind of observation was one of the crucial ingredients in Soundararajan's work~\cite{soundmoments}: see Proposition 1, below.
\end{enumerate}

When thinking about the Dirichlet polynomial contribution $\Re \sum_{p \leq T} \frac{1}{p^{1/2+it}}$ as $T \leq t \leq 2T$ varies, a standard heursitic is that its average behaviour should be the same as the behaviour of
$$ \Re \sum_{p \leq T} \frac{U_{p}}{p^{1/2}} , $$
where $U_{p}$ are independent random variables distributed uniformly on the unit circle in $\C$. The latter object is a real valued random variable with mean zero and variance $(1/2) \sum_{p \leq T} 1/p \approx (1/2)\log\log T$. As Soundararajan~\cite{soundmoments} discusses, if we knew that $\log|\zeta(1/2+it)|$ behaved in the same way as a Gaussian random variable with such mean and variance then we could obtain a moment bound as in Theorem 1, and his argument can be seen as an attempt to demonstrate as much Gaussian-like behaviour as possible for $\Re \sum_{p \leq T} \frac{1}{p^{1/2+it}}$, thereby coming close to a sharp moment bound. We will push this line of thought further and obtain the sharp bound in Theorem 1. Before explaining in detail, we note something else about the random object $\Re \sum_{p \leq T} \frac{U_{p}}{p^{1/2}} $ : given any values $1 = x_{0} < x_{1} < x_{2} < ... < x_{n} = T$ we can split
$$ \Re \sum_{p \leq T} \frac{U_{p}}{p^{1/2}} = \sum_{i=1}^{n} \Re \sum_{x_{i-1} < p \leq x_{i}} \frac{U_{p}}{p^{1/2}} , $$
where the pieces $\Re \sum_{x_{i-1} < p \leq x_{i}} \frac{U_{p}}{p^{1/2}}$ are independent of one another, and have mean zero and variance $(1/2) \sum_{x_{i-1} < p \leq x_{i}} 1/p \approx (1/2) \log(\log x_{i} /\log x_{i-1})$. In particular, the latter terms in the sum may not contribute much to the variance, and so will not contribute much to the typical size.

\vspace{12pt}
Soundararajan's argument~\cite{soundmoments} actually works by upper bounding $\log|\zeta(1/2+it)|$ by a Dirichlet polynomial $A(1/2+it)$ of suitably chosen length, and investigating the frequency of large values of such Dirichlet polynomials by Markov-inequality type arguments applied to their high moments:
$$ \text{meas}\{T \leq t \leq 2T : |A(1/2+it)| \geq V\} \leq \frac{1}{V^{U}} \int_{T}^{2T} |A(1/2+it)|^{U} dt , \;\;\;\;\; \forall V, U \geq 0 . $$
There are two reasons why this doesn't yield sharp estimates:
\begin{enumerate}
\item If one is studying a Dirichlet polynomial $A(s) = \sum_{1 \leq n \leq T^{1/W}} \frac{a_{n}}{n^{s}}$ of length $T^{1/W}$, one can only obtain good bounds for the first $O(W)$ moments. So to study high moments one must work with short Dirichlet polynomials, which produces an error term in the upper bound for $\log|\zeta(1/2+it)|$. See Proposition 1, below.

\item One cannot recover optimal bounds for the frequency of large values just by Markov inequality-type arguments using moments\footnote{For example, the reader may check that if we assume Theorem 1 for all $k$, and if $\log\log T \leq V \leq 1000\log\log T$ (say), then the best bound we can obtain on $\text{meas}\{T \leq t \leq 2T : \log|\zeta(1/2+it)| \geq V\}$ using Markov's inequality is $\ll T e^{-V^{2}/\log\log T}$ (by choosing $k = V/\log\log T$). The true answer on this range is presumably $\ll (T/\sqrt{\log\log T}) e^{-V^{2}/\log\log T}$.}.
\end{enumerate}

Roughly speaking, to address the first problem we will upper bound $\log|\zeta(1/2+it)|$ by a long Dirichlet polynomial, but split this polynomial into multiple pieces and raise each piece to a different power (and then look at the integral of the product of all the powers of polynomials). As discussed previously, in the upper bound for $\log|\zeta(1/2+it)|$ one expects later terms in the Dirichlet polynomial to contribute increasingly little to the total size, so one can afford to raise the later (longer) pieces only to the smaller powers that are permitted by their increased length\footnote{An idea like this already occurs in Soundararajan's work~\cite{soundmoments}, when he splits his Dirichlet polynomial $A(s)$ into two pieces. But Soundararajan examines each piece separately, whereas we study the joint behaviour of all of our polynomials.}. This plan can only be implemented if one knows in advance that the various pieces of the Dirichlet polynomial take roughly their expected size, but if $t$ is a ``bad'' point for which a piece is too big then one can use a truncated Dirichlet polynomial upper bound for $\log|\zeta(1/2+it)|$, and still win because the set of such $t$ has small measure.

To address the second problem, we will not actually estimate the frequency of large values of $\log|\zeta(1/2+it)|$ at all, but will work directly with moment-type objects throughout. Thus we will have, roughly speaking,
$$ |\zeta(1/2+it)|^{2k} = \exp(2k\log|\zeta(1/2+it)|) \leq \exp(2k\Re A(1/2+it)) = \prod_{j} \exp(2k\Re A_{j}(1/2+it)) , $$
where $A_{j}(s)$ are Dirichlet polynomials that sum to give $A(s)$. If the size of each piece $A_{j}(1/2+it)$ is well controlled then one can replace each factor $\exp(2k\Re A_{j}(1/2+it))$ by a sum of a few terms from its series expansion, and then bound the integral of the resulting product, which becomes a sum of products of Dirichlet polynomials. This last manoeuvre was inspired by a very nice paper of Radziwi\l\l~\cite{rad2} on the Selberg central limit theorem, in which it is used to estimate the moment generating function of (something like) $\log|\zeta(1/2+it)|$.

\vspace{12pt}
As in the proof of the conditional estimate $\int_{T}^{2T} |\zeta(1/2+it)|^{2k} dt \ll_{k, \epsilon} T \log^{k^{2}+\epsilon}T $, our proof of Theorem 1 is very general and extends to give (presumably) sharp upper bounds\footnote{That is, our method gives bounds that match the conjectured behaviour. We do not yet have matching lower bounds in all cases when $k$ is small: see e.g. the papers of Chandee and Li~\cite{chandeeli}, Radziwi\l{\l} and Soundararajan~\cite{radsound}, and Rudnick and Soundararajan~\cite{rudnicksound} for some state of the art lower bound results.} for other moments of $L$-functions. We discuss this in $\S 5$, below. Straightforward adaptations would also yield sharp bounds for the short interval moments $\int_{T}^{T + T^{\theta}} |\zeta(1/2+it)|^{2k} dt$ considered by Ivi\'{c}~\cite{ivicintervals}, where $0 < \theta \leq 1$ and $k \geq 0$.

The last section of this paper is a discussion about the implicit constant in Theorem 1. Our proof supplies an implicit constant of the form $e^{e^{O(k)}}$, for large $k$ and $T$, whereas random matrix theory predicts a constant of size around $e^{-k^{2}\log k}$. In $\S 6$ we explain how, conjecturally, one could refine our method to obtain an implicit constant of size about $e^{-k^{2}\log k}$, {\em without} using random matrix theory.

\vspace{12pt}
Finally, we remark that in some work in progress, Radziwi\l{\l} and Soundararajan have independently developed a splitting technique similar to the one used here. They will apply this to prove sharp lower bounds for moments of $L$-functions, as well as unconditional sharp upper bounds for the moments of $\zeta(s)$ with $0 \leq k \leq 2$ real.

\section{Some tools}
In this section we recall two facts that will be fundamental tools for the proof of Theorem 1.

We begin with the following result, that gives (conditional on the Riemann Hypothesis) an upper bound for $\log|\zeta(1/2+it)|$ in terms of a Dirichlet polynomial. In fact the result allows considerable freedom in choosing the length of that Dirichlet polynomial, and that will be very important for our arguments.
\begin{prop1}[Adapted from the Main Proposition of Soundararajan~\cite{soundmoments}]
Assume the Riemann Hypothesis is true, and let $T$ be large. Then for any $2 \leq x \leq T^{2}$, and any $T \leq t \leq 2T$, we have
$$ \log|\zeta(1/2+it)| \leq \Re\left(\sum_{p \leq x} \frac{1}{p^{1/2+1/\log x + it}} \frac{\log(x/p)}{\log x} + \sum_{p \leq \min\{\sqrt{x},\log T\}} \frac{(1/2)}{p^{1 + 2it}} \right) + \frac{\log T}{\log x} + O(1) , $$
where $p$ denotes primes.
\end{prop1}

Proposition 1 follows by choosing $\lambda = 1$ in the Main Proposition of Soundararajan~\cite{soundmoments}, and noting that the contribution from prime cubes and higher powers there is $O(1)$, and the contribution from prime squares there is
$$ \Re \sum_{p \leq \sqrt{x}} \frac{(1/2)}{p^{1 + 2/\log x + 2it}} \frac{\log(x/p^{2})}{\log x} = \Re \sum_{p \leq \sqrt{x}} \frac{(1/2)}{p^{1 + 2it}} + O(1) = \Re \sum_{p \leq \min\{\sqrt{x},\log T\}} \frac{(1/2)}{p^{1 + 2it}} + O(1) . $$
Here the truncation of the sum at $\log T$ is easily justified, assuming the Riemann Hypothesis, by standard explicit formula arguments (truncating first at $\log^{10}T$, and then observing that $\sum_{\log T < p \leq \log^{10}T} 1/p = O(1)$): see e.g. the proof of Lemma 2 of Soundararajan~\cite{soundmoments}.

We will also need quite precise information about the integral of a product of terms $\cos(t\log p)$. Our formulation of this is a slight variant of Lemma 4 from Radziwi\l\l's paper~\cite{rad2} on Selberg's central limit theorem, noting that the error term there is estimated quite generously, and can actually be taken as $O(n)$.
\begin{prop2}
Let $T$ be large and let $n=p_{1}^{\alpha_{1}}...p_{r}^{\alpha_{r}}$, where the $p_{i}$ are distinct primes and $\alpha_{i} \in \N$ for all $i$. Then
$$ \int_{T}^{2T} \prod_{i=1}^{r} (\cos(t\log p_{i}))^{\alpha_{i}} dt = Tf(n) + O(n) , $$
where $f(n)=0$ if any of the exponents $\alpha_{i}$ is odd, and otherwise
$$ f(n) := \prod_{i=1}^{r} \frac{1}{2^{\alpha_{i}}} \frac{\alpha_{i}!}{((\alpha_{i}/2)!)^{2}} . $$
\end{prop2}
The fact that there is no main term here unless all the exponents are even, meaning that the primes match in pairs, may be familiar to the reader from other moment computations, and provides encouragement that we will see nice (e.g. Gaussian-like) behaviour in our calculations. The fact that $f(n)$ is a multiplicative function, so that distinct primes do not interact with each other, is also encouraging, since it reflects our expectation that the quantities $\cos(t\log p) = \Re p^{it}$ behave ``independently'', as $t$ varies, for distinct primes.

\section{Outline of the proof of Theorem 1}
Suppose, as always, that $T$ is large, and that $k \geq 1$ is a fixed real number. We could also handle the case of $0 \leq k < 1$ by an easy adaptation of our argument, but since Theorem 1 is already known on that range we will ignore it, because it will simplify matters if we can assume that multiplying by $k$ doesn't make terms smaller. In this section we will set up some basic notation, state three lemmas, and deduce Theorem 1 from those lemmas. The lemmas will be proved in the next section.

\vspace{12pt}
Firstly, define the sequence $(\beta_{i})_{i \geq 0}$ by
$$ \beta_{0} := 0, \;\;\;\;\; \beta_{i} := \frac{20^{i-1}}{(\log\log T)^{2}} \;\;\; \forall i \geq 1, $$
and define
$$ \mathcal{I} = \mathcal{I}_{k,T} := 1 + \max\{i : \beta_{i} \leq e^{-1000k}\} . $$

Next, for the sake of concision we shall introduce notation for certain Dirichlet polynomials. For each $1 \leq i \leq j \leq \mathcal{I}$, set
$$ G_{(i,j)}(t) = G_{(i,j),T}(t) := \sum_{T^{\beta_{i-1}} < p \leq T^{\beta_{i}}} \frac{1}{p^{1/2 + 1/(\beta_{j} \log T) + it}} \frac{\log(T^{\beta_{j}}/p)}{\log(T^{\beta_{j}})} . $$

Finally, let us define the set $\mathcal{T} = \mathcal{T}_{k,T}$, by
$$ \mathcal{T} := \{T \leq t \leq 2T : \left| \Re \sum_{T^{\beta_{i-1}} < p \leq T^{\beta_{i}}} \frac{1}{p^{1/2 + 1/(\beta_{\mathcal{I}} \log T) + it}} \frac{\log(T^{\beta_{\mathcal{I}}}/p)}{\log(T^{\beta_{\mathcal{I}}})} \right| \leq \beta_{i}^{-3/4} \; \forall 1 \leq i \leq \mathcal{I} \}, $$
and let us define sets $\mathcal{S}(j) = \mathcal{S}_{k,T}(j)$, for $0 \leq j \leq \mathcal{I} - 1$, by
$$ \mathcal{S}(j) := \{T \leq t \leq 2T : |\Re G_{(i,l)}(t)| \leq \beta_{i}^{-3/4} \; \forall 1 \leq i \leq j, \; \forall i \leq l \leq \mathcal{I}, $$
$$\;\;\;\;\; \text{but } |\Re G_{(j+1,l)}(t)| > \beta_{j+1}^{-3/4} \; \text{ for some } j+1 \leq l \leq \mathcal{I} \} . $$

\begin{lem1}
Let the situation be as described above, where $k \geq 1$ and $T \geq e^{e^{(10000k)^{2}}}$ is large enough. Then
$$ \int_{t \in \mathcal{T}} \exp\left( 2k \Re \sum_{p \leq T^{\beta_{\mathcal{I}}}} \frac{1}{p^{1/2 + 1/(\beta_{\mathcal{I}} \log T) + it}} \frac{\log(T^{\beta_{\mathcal{I}}}/p)}{\log(T^{\beta_{\mathcal{I}}})} \right) dt \ll T \log^{k^{2}}T , $$
where the implicit constant is absolute.
\end{lem1}

\begin{lem2}
Let the situation be as described above, where $k \geq 1$ and $T \geq e^{e^{(10000k)^{2}}}$ is large enough. Then $\text{meas}(\mathcal{S}(0)) \ll T e^{-(\log\log T)^{2}/10}$, and for any $1 \leq j \leq \mathcal{I} - 1$ we have
$$ \int_{t \in \mathcal{S}(j)} \exp\left( 2k \Re \sum_{p \leq T^{\beta_{j}}} \frac{1}{p^{1/2 + 1/(\beta_{j} \log T) + it}} \frac{\log(T^{\beta_{j}}/p)}{\log(T^{\beta_{j}})} \right) dt \ll e^{-\beta_{j+1}^{-1} \log(1/\beta_{j+1}) /21} T \log^{k^{2}}T , $$
where the implicit constant is absolute.
\end{lem2}

\begin{lem3}
The estimates in Lemma 1, and in Lemma 2 for $1 \leq j \leq \mathcal{I}-1$, remain true when the Dirichlet polynomials in the exponents are replaced by
$$ \sum_{p \leq T^{\beta_{j}}} \frac{1}{p^{1/2 + 1/(\beta_{j} \log T) + it}} \frac{\log(T^{\beta_{j}}/p)}{\log(T^{\beta_{j}})} + \sum_{p \leq \log T} \frac{(1/2)}{p^{1 + 2it}} , $$
except that the implicit constants in the estimates may now depend on $k$.
\end{lem3}

Now we can swiftly deduce Theorem 1. We note that
$$ [T,2T] = \bigcup_{j=0}^{\mathcal{I}-1} \mathcal{S}(j) \cup \mathcal{T} , $$
so it will suffice to show that
$$ \sum_{j=0}^{\mathcal{I}-1} \int_{t \in \mathcal{S}(j)} |\zeta(1/2+it)|^{2k} dt + \int_{t \in \mathcal{T}} |\zeta(1/2+it)|^{2k} dt \ll_{k} T \log^{k^{2}}T . $$
But Proposition 1 implies that
$$ \log|\zeta(1/2+it)| \leq \Re \left(\sum_{p \leq T^{\beta_{\mathcal{I}}}} \frac{1}{p^{1/2 + 1/(\beta_{\mathcal{I}} \log T) + it}} \frac{\log(T^{\beta_{\mathcal{I}}}/p)}{\log(T^{\beta_{\mathcal{I}}})} + \sum_{p \leq \log T} \frac{1/2}{p^{1 + 2it}} \right) + \frac{1}{\beta_{\mathcal{I}}} + O(1) , $$
and therefore $ \int_{t \in \mathcal{T}} |\zeta(1/2+it)|^{2k} dt$ is
$$ \ll_{k} e^{2k/\beta_{\mathcal{I}}} \int_{t \in \mathcal{T}} \exp\left( 2k \Re \left(\sum_{p \leq T^{\beta_{\mathcal{I}}}} \frac{1}{p^{1/2 + 1/(\beta_{\mathcal{I}} \log T) + it}} \frac{\log(T^{\beta_{\mathcal{I}}}/p)}{\log(T^{\beta_{\mathcal{I}}})} + \sum_{p \leq \log T} \frac{1/2}{p^{1 + 2it}} \right)  \right) dt \ll_{k} T \log^{k^{2}}T , $$
using Lemma 3. Similarly, if $1 \leq j \leq \mathcal{I}-1$ then Proposition 1 implies that
$$ \log|\zeta(1/2+it)| \leq \Re \left(\sum_{p \leq T^{\beta_{j}}} \frac{1}{p^{1/2 + 1/(\beta_{j} \log T) + it}} \frac{\log(T^{\beta_{j}}/p)}{\log(T^{\beta_{j}})} + \sum_{p \leq \log T} \frac{1/2}{p^{1 + 2it}} \right) + \frac{1}{\beta_{j}} + O(1) , $$
and so Lemma 3 implies that
$$ \int_{t \in \mathcal{S}(j)} |\zeta(1/2+it)|^{2k} dt \ll_{k} e^{2k/\beta_{j}} \cdot e^{-\beta_{j+1}^{-1} \log(1/\beta_{j+1}) /21} T \log^{k^{2}}T . $$
And we see
\begin{equation}\label{sumlemma2}
e^{2k/\beta_{j}} \cdot e^{-\beta_{j+1}^{-1} \log(1/\beta_{j+1}) /21} = e^{2k/\beta_{j} - \log(1/\beta_{j+1}) /(420\beta_{j})} \leq e^{-0.01k/\beta_{j}},
\end{equation}
since $\beta_{j+1} \leq \beta_{\mathcal{I}} \leq 20 e^{-1000k}$, which implies that $\log(1/\beta_{j+1}) \geq 900k$ (say). Then the sum of these bounds over $1 \leq j \leq \mathcal{I}-1$ is also $ \ll_{k} T \log^{k^{2}}T$. 

Finally, when $j=0$ we note that
\begin{eqnarray}
\int_{t \in \mathcal{S}(0)} |\zeta(1/2+it)|^{2k} dt & \leq & \sqrt{\text{meas}(\mathcal{S}(0)) \cdot \int_{T}^{2T} |\zeta(1/2+it)|^{4k} dt} \nonumber \\
& \ll_{k} & \sqrt{T e^{-(\log\log T)^{2}/10} \cdot T \log^{(2k)^{2}+1}T} , \nonumber
\end{eqnarray}
in view of Lemma 2 and Soundararajan's bound for the moments of the zeta function (with $\epsilon = 1$). This is certainly $\ll_{k} T\log^{k^{2}}T$, and so our proof of Theorem 1 is complete.
\begin{flushright}
Q.E.D.
\end{flushright}

\section{Proofs of the Lemmas}

\subsection{Proof of Lemma 1}
Our goal is to show that
$$ \int_{t \in \mathcal{T}} \exp\left( 2k \Re \sum_{p \leq T^{\beta_{\mathcal{I}}}} \frac{1}{p^{1/2 + 1/(\beta_{\mathcal{I}} \log T) + it}} \frac{\log(T^{\beta_{\mathcal{I}}}/p)}{\log(T^{\beta_{\mathcal{I}}})} \right) dt \ll T \log^{k^{2}}T . $$
For the sake of concision we shall set
$$ F_{i}(t) = F_{i,k,T}(t) := \sum_{T^{\beta_{i-1}} < p \leq T^{\beta_{i}}} \frac{1}{p^{1/2 + 1/(\beta_{\mathcal{I}} \log T) + it}} \frac{\log(T^{\beta_{\mathcal{I}}}/p)}{\log(T^{\beta_{\mathcal{I}}})} \;\;\; \forall 1 \leq i \leq \mathcal{I} , $$
(note this is just an alternative notation for $G_{(i,\mathcal{I})}(t)$), so the integral we are trying to bound is
$$  \int_{t \in \mathcal{T}} \prod_{1 \leq i \leq \mathcal{I}} \exp\left( 2k \Re F_{i}(t) \right) dt = \int_{t \in \mathcal{T}} \prod_{1 \leq i \leq \mathcal{I}} \left( \exp\left( k \Re F_{i}(t) \right) \right)^{2} dt  . $$
Using the series expansion of the exponential function, and recalling that $|\Re F_{i}(t)| \leq \beta_{i}^{-3/4}$ when $t \in \mathcal{T}$, which means that we can truncate the series with a very small error, we find
\begin{eqnarray}
&& \int_{t \in \mathcal{T}} \prod_{1 \leq i \leq \mathcal{I}} \exp\left( 2k \Re F_{i}(t) \right) dt \nonumber \\
& = & \int_{t \in \mathcal{T}} \prod_{1 \leq i \leq \mathcal{I}} \left( \sum_{0 \leq j \leq 100k \beta_{i}^{-3/4}} \frac{(k \Re F_{i}(t))^{j}}{j!} + O\left(\frac{(k |\Re F_{i}(t)|)^{[100k\beta_{i}^{-3/4}] + 1}}{([100k\beta_{i}^{-3/4}] + 1)!}\right) \right)^{2} dt  \nonumber \\
& = & \int_{t \in \mathcal{T}} \prod_{1 \leq i \leq \mathcal{I}} \left(1+ O(e^{-100k\beta_{i}^{-3/4}}) \right) \left( \sum_{0 \leq j \leq 100k \beta_{i}^{-3/4}} \frac{(k \Re F_{i}(t))^{j}}{j!} \right)^{2} dt \nonumber \\
& \leq & \left(1+ O(e^{-100k\beta_{\mathcal{I}}^{-3/4}}) \right) \int_{T}^{2T} \prod_{1 \leq i \leq \mathcal{I}} \left( \sum_{0 \leq j \leq 100k \beta_{i}^{-3/4}} \frac{(k \Re F_{i}(t))^{j}}{j!} \right)^{2} dt .  \nonumber
\end{eqnarray}
Here $[\cdot]$ denotes the integer part of $\cdot$, and in the final line we completed the range of integration to the entire interval $[T,2T]$, which is permissible since we are seeking an upper bound and the integrand is always positive. Obtaining this positivity is the reason that we wrote $\exp\left( 2k \Re F_{i}(t) \right) = \left( \exp\left( k \Re F_{i}(t) \right) \right)^{2}$ in the first place, rather than looking directly at the series expansion of $\exp\left( 2k \Re F_{i}(t) \right)$.

Now if we expand all of the $j$-th powers, the squares, and the product over $i$, and recall that
$$ \Re F_{i}(t) = \sum_{T^{\beta_{i-1}} < p \leq T^{\beta_{i}}} \frac{\cos(t\log p)}{p^{1/2 + 1/(\beta_{\mathcal{I}} \log T)}} \frac{\log(T^{\beta_{\mathcal{I}}}/p)}{\log(T^{\beta_{\mathcal{I}}})} \;\;\; \forall 1 \leq i \leq \mathcal{I}, \;\;\; T \leq t \leq 2T, $$
we see the integral above is equal to
$$ \sum_{\tilde{j}, \tilde{l}} \left( \prod_{1 \leq i \leq \mathcal{I}} \frac{k^{j_{i}}}{j_{i}!} \frac{k^{l_{i}}}{l_{i}!} \right) \sum_{\tilde{p}, \tilde{q}} C(\tilde{p},\tilde{q}) \int_{T}^{2T} \prod_{1 \leq i \leq \mathcal{I}} \left( \prod_{\substack{1 \leq r \leq j_{i}, \\ 1 \leq s \leq l_{i}}} \cos(t\log p(i,r)) \cos(t\log q(i,s)) \right) dt , $$
where the outer sum is over all vectors $\tilde{j} = (j_{1},j_{2},...,j_{\mathcal{I}}), \tilde{l} = (l_{1},l_{2},...,l_{\mathcal{I}})$ with components satisfying
$$ 0 \leq j_{i}, l_{i} \leq 100k\beta_{i}^{-3/4}; $$
the inner sum is over all vectors $\tilde{p} = (p(1,1),p(1,2),...,p(1,j_{1}),p(2,1),...,p(2,j_{2}),...,p(\mathcal{I},j_{\mathcal{I}}))$, $\tilde{q} = (q(1,1),...,q(\mathcal{I},l_{\mathcal{I}}))$ with components that are primes satisfying
$$ T^{\beta_{i-1}} < p(i,1),...,p(i,j_{i}), q(i,1),...,q(i,l_{i}) \leq T^{\beta_{i}} \;\;\; \forall 1 \leq i \leq \mathcal{I} ; $$
and
$$ C(\tilde{p},\tilde{q}) := \prod_{1 \leq i \leq \mathcal{I}} \left( \prod_{\substack{1 \leq r \leq j_{i}, \\ 1 \leq s \leq l_{i}}} \frac{1}{p(i,r)^{1/2 + 1/(\beta_{\mathcal{I}} \log T)}} \frac{\log(T^{\beta_{\mathcal{I}}}/p(i,r))}{\log(T^{\beta_{\mathcal{I}}})} \frac{1}{q(i,s)^{1/2 + 1/(\beta_{\mathcal{I}} \log T)}} \frac{\log(T^{\beta_{\mathcal{I}}}/q(i,s))}{\log(T^{\beta_{\mathcal{I}}})} \right) . $$
At this point we can note that in the above setting,
$$ \prod_{1 \leq i \leq \mathcal{I}} \prod_{\substack{1 \leq r \leq j_{i}, \\ 1 \leq s \leq l_{i}}} p(i,r) q(i,s) \leq \prod_{1 \leq i \leq \mathcal{I}} T^{\beta_{i}(j_{i}+l_{i})} \leq \prod_{1 \leq i \leq \mathcal{I}} T^{200k\beta_{i}^{1/4}} \leq T^{400k\beta_{\mathcal{I}}^{1/4}} \leq T^{0.1}, $$
say, since the numbers $\beta_{i}^{1/4}$ form a geometric progression with common ratio $20^{1/4} \geq 2$, and $\beta_{\mathcal{I}}^{1/4} \leq \beta_{1}^{1/4} + (20e^{-1000k})^{1/4} = (\log\log T)^{-1/2} + (20e^{-1000k})^{1/4} \leq 1/(5000k)$ (since we assume that $T \geq e^{e^{(10000k)^{2}}}$). In view of Proposition 2, this means that
$$ \int_{T}^{2T} \prod_{1 \leq i \leq \mathcal{I}} \left( \prod_{\substack{1 \leq r \leq j_{i}, \\ 1 \leq s \leq l_{i}}} \cos(t\log p(i,r)) \cos(t\log q(i,s)) \right) dt = T f(\prod_{1 \leq i \leq \mathcal{I}} \prod_{\substack{1 \leq r \leq j_{i}, \\ 1 \leq s \leq l_{i}}} p(i,r) q(i,s)) + O(T^{0.1}) . $$
Since $f(n)$ is a non-negative function, and $C(\tilde{p},\tilde{q})$ is at most as large as the rather simpler quantity $ D(\tilde{p},\tilde{q}) := \prod_{1 \leq i \leq \mathcal{I}} \left( \prod_{\substack{1 \leq r \leq j_{i}, \\ 1 \leq s \leq l_{i}}} \frac{1}{\sqrt{p(i,r)}} \frac{1}{\sqrt{q(i,s)}} \right) $, we can deduce that
\begin{eqnarray}
\int_{t \in \mathcal{T}} \prod_{1 \leq i \leq \mathcal{I}} \exp\left( 2k \Re F_{i}(t) \right) dt & \ll & T \sum_{\tilde{j}, \tilde{l}} \left( \prod_{1 \leq i \leq \mathcal{I}} \frac{k^{j_{i}}}{j_{i}!} \frac{k^{l_{i}}}{l_{i}!} \right) \sum_{\tilde{p}, \tilde{q}} D(\tilde{p},\tilde{q}) f(\prod_{1 \leq i \leq \mathcal{I}} \prod_{\substack{1 \leq r \leq j_{i}, \\ 1 \leq s \leq l_{i}}} p(i,r) q(i,s)) \nonumber \\
&& + T^{0.1} \sum_{\tilde{j}, \tilde{l}} \left( \prod_{1 \leq i \leq \mathcal{I}} \frac{k^{j_{i}}}{j_{i}!} \frac{k^{l_{i}}}{l_{i}!} \right) \sum_{\tilde{p}, \tilde{q}} D(\tilde{p},\tilde{q}) . \nonumber
\end{eqnarray}
The second term here can be rewritten as
\begin{eqnarray}
T^{0.1} \prod_{1 \leq i \leq \mathcal{I}} \left( \sum_{0 \leq j \leq 100k\beta_{i}^{-3/4}} \frac{k^{j}}{j!} \left(\sum_{T^{\beta_{i-1}} < p \leq T^{\beta_{i}}} \frac{1}{\sqrt{p}} \right)^{j} \right)^{2} & \leq & T^{0.1} \prod_{1 \leq i \leq \mathcal{I}} T^{200k\beta_{i}^{1/4}} \left( \sum_{0 \leq j \leq 100k\beta_{i}^{-3/4}} \frac{k^{j}}{j!} \right)^{2} \nonumber \\
& \leq & T^{0.2} e^{2k\mathcal{I}}, \nonumber
\end{eqnarray}
and this is certainly $\leq T^{0.2} (\log\log T)^{2k}$, which is negligible.

Finally, since $f$ is a multiplicative (though not totally multiplicative) function we can now reassemble the foregoing (rather horrible!) bound into a product. Indeed, we find that
\begin{eqnarray}\label{neatbound}
&& \int_{t \in \mathcal{T}} \prod_{1 \leq i \leq \mathcal{I}} \exp\left( 2k \Re F_{i}(t) \right) dt \\
& \ll & T \prod_{1 \leq i \leq \mathcal{I}} \sum_{0 \leq j,l \leq 100k\beta_{i}^{-3/4}} \frac{k^{j + l}}{j! l!} \sum_{T^{\beta_{i-1}} < p_{1},...,p_{j}, q_{1},...,q_{l} \leq T^{\beta_{i}}} \frac{f(p_{1}...p_{j}q_{1}...q_{l})}{\sqrt{p_{1}...p_{j}q_{1}...q_{l}}} + T^{0.2} (\log\log T)^{2k}, \nonumber
\end{eqnarray}
and the first term here is clearly equal to
$$ T \prod_{1 \leq i \leq \mathcal{I}} \sum_{0 \leq m \leq 200k\beta_{i}^{-3/4}} k^{m} \left(\sum_{\substack{j+l=m, \\ 0 \leq j,l \leq 100k\beta_{i}^{-3/4}}} \frac{1}{j! l!} \right) \sum_{T^{\beta_{i-1}} < p_{1},...,p_{m} \leq T^{\beta_{i}}} \frac{f(p_{1}...p_{m})}{\sqrt{p_{1}...p_{m}}}, $$
which is
$$ \leq T \prod_{1 \leq i \leq \mathcal{I}} \sum_{0 \leq m \leq 200k\beta_{i}^{-3/4}} \frac{k^{m} 2^{m}}{m!} \sum_{T^{\beta_{i-1}} < p_{1},...,p_{m} \leq T^{\beta_{i}}} \frac{f(p_{1}...p_{m})}{\sqrt{p_{1}...p_{m}}} . $$
The reader may note that we have now undone the ``introduction of squares'' (by splitting $\exp\left( 2k \Re F_{i}(t) \right)$ as $\left( \exp\left( k \Re F_{i}(t) \right) \right)^{2}$) that we performed earlier, and recovered an expression involving powers of $2k$.

The function $f$ is supported on squares, so we can certainly restrict the sum over $m$ to even terms $m=2n$. Thus the above expression is equal to
$$ T \prod_{1 \leq i \leq \mathcal{I}} \sum_{0 \leq n \leq 100k\beta_{i}^{-3/4}} \frac{(2k)^{2n}}{(2n)!} \sum_{T^{\beta_{i-1}} < p_{1},...,p_{n} \leq T^{\beta_{i}}} \frac{f(p_{1}^{2}...p_{n}^{2})}{p_{1}...p_{n}} \frac{\#\{(q_{1},...,q_{2n}) : q_{1}...q_{2n} = p_{1}^{2}...p_{n}^{2} \}}{\#\{(q_{1},...,q_{n}) : q_{1}...q_{n} = p_{1}...p_{n} \}}, $$
where $q_{1},...,q_{2n}$ again denote primes between $T^{\beta_{i-1}}$ and $T^{\beta_{i}}$. (Here the factor involving the $q_{i}$ is a correction factor that ensures we count each integer $p_{1}^{2}...p_{n}^{2}$ the correct number of times.) Now if $p_{1}...p_{n}$ is a product of $r$ distinct primes, with multiplicities $\alpha_{1},...,\alpha_{r}$ (so that $\alpha_{1} + ... + \alpha_{r} = n$), then we have
\begin{equation}
f(p_{1}^{2}...p_{n}^{2}) = \frac{1}{2^{2n}} \prod_{j=1}^{r} \frac{(2\alpha_{j})!}{(\alpha_{j}!)^{2}}, \;\;\; \textrm{ and } \;\;\; \#\{(q_{1},...,q_{2n}) : q_{1}...q_{2n} = p_{1}^{2}...p_{n}^{2} \} = \frac{(2n)!}{\prod_{j=1}^{r} (2\alpha_{j})!}, \nonumber
\end{equation}
\begin{equation}\label{countingtuples}
\textrm{and } \;\;\; \#\{(q_{1},...,q_{n}) : q_{1}...q_{n} = p_{1}...p_{n} \} = \frac{n!}{\prod_{j=1}^{r} \alpha_{j}!} ,
\end{equation}
and so the above expression is equal to
\begin{eqnarray}
&& T \prod_{1 \leq i \leq \mathcal{I}} \sum_{0 \leq n \leq 100k\beta_{i}^{-3/4}} \frac{k^{2n}}{n!} \sum_{T^{\beta_{i-1}} < p_{1},...,p_{n} \leq T^{\beta_{i}}} \frac{1}{p_{1}...p_{n}} \frac{1}{\prod_{j=1}^{r} \alpha_{j}!} \nonumber \\
& \leq & T \prod_{1 \leq i \leq \mathcal{I}} \sum_{0 \leq n \leq 100k\beta_{i}^{-3/4}} \frac{1}{n!} \left(k^{2} \sum_{T^{\beta_{i-1}} < p \leq T^{\beta_{i}}} \frac{1}{p} \right)^{n} \leq T \exp\left(k^{2} \sum_{p \leq T^{\beta_{\mathcal{I}}}} \frac{1}{p} \right) . \nonumber
\end{eqnarray}
This is indeed $\leq T \exp(k^{2}\log\log T) = T\log^{k^{2}}T$, as claimed in Lemma 1.
\begin{flushright}
Q.E.D.
\end{flushright}

\subsection{Proof of Lemma 2}
The proof of Lemma 2 is a direct modification of the proof of Lemma 1, so we shall simply outline the main details. We have
\begin{eqnarray}
&& \int_{t \in \mathcal{S}(j)} \exp\left( 2k \Re \sum_{p \leq T^{\beta_{j}}} \frac{1}{p^{1/2 + 1/(\beta_{j} \log T) + it}} \frac{\log(T^{\beta_{j}}/p)}{\log(T^{\beta_{j}})} \right) dt \nonumber \\
& \leq & \sum_{l=j+1}^{\mathcal{I}} \int_{\substack{T \leq t \leq 2T: |\Re G_{(i,j)}(t)| \leq \beta_{i}^{-3/4} \; \forall 1 \leq i \leq j, \\ \text{but } |\Re G_{(j+1,l)}(t)| > \beta_{j+1}^{-3/4}}} \exp\left( 2k \Re \sum_{p \leq T^{\beta_{j}}} \frac{1}{p^{1/2 + 1/(\beta_{j} \log T) + it}} \frac{\log(T^{\beta_{j}}/p)}{\log(T^{\beta_{j}})} \right) dt  \nonumber \\
& = & \sum_{l=j+1}^{\mathcal{I}} \int_{\substack{T \leq t \leq 2T: |\Re G_{(i,j)}(t)| \leq \beta_{i}^{-3/4} \; \forall 1 \leq i \leq j, \\ \text{but } |\Re G_{(j+1,l)}(t)| > \beta_{j+1}^{-3/4}}} \prod_{1 \leq i \leq j} \left( \exp\left( k \Re G_{(i,j)}(t) \right) \right)^{2} dt  \nonumber ,
\end{eqnarray}
by definition of the set $\mathcal{S}(j)$; and each of the integrals here is
\begin{eqnarray}
& \leq & \int_{T \leq t \leq 2T: |\Re G_{(i,j)}(t)| \leq \beta_{i}^{-3/4} \; \forall 1 \leq i \leq j} \prod_{1 \leq i \leq j} \left( \exp\left( k \Re G_{(i,j)}(t) \right) \right)^{2}  \left(\beta_{j+1}^{3/4} \Re G_{(j+1,l)}(t) \right)^{2[1/(10\beta_{j+1})]}  dt \nonumber \\
& \ll & (\beta_{j+1}^{3/2})^{[1/(10\beta_{j+1})]} \int_{T}^{2T} \prod_{1 \leq i \leq j} \left( \sum_{0 \leq n \leq 100k \beta_{i}^{-3/4}} \frac{(k \Re G_{(i,j)}(t))^{n}}{n!} \right)^{2} \left( \Re G_{(j+1,l)}(t) \right)^{2[1/(10\beta_{j+1})]}  dt , \nonumber
\end{eqnarray}
as in the proof of Lemma 1. For the sake of concision, let us temporarily set $M := 2[1/(10\beta_{j+1})]$.

Now expanding the squares and the powers of $n$, and proceeding as in the proof of Lemma 1 (see e.g. lines (\ref{neatbound}) and (\ref{countingtuples}), and the surrounding discussion), we find that
\begin{eqnarray}
&& \int_{T}^{2T} \prod_{1 \leq i \leq j} \left( \sum_{0 \leq n \leq 100k \beta_{i}^{-3/4}} \frac{(k \Re G_{(i,j)}(t))^{n}}{n!} \right)^{2} \left( \Re G_{(j+1,l)}(t) \right)^{2[1/(10\beta_{j+1})]}  dt \nonumber \\
& \ll & T \prod_{1 \leq i \leq j} \left( \sum_{0 \leq m \leq 200k\beta_{i}^{-3/4}} \frac{k^{m} 2^{m}}{m!} \sum_{T^{\beta_{i-1}} < p_{1},...,p_{m} \leq T^{\beta_{i}}} \frac{f(p_{1}...p_{m})}{\sqrt{p_{1}...p_{m}}} \right) \cdot \sum_{T^{\beta_{j}} < p_{1},...,p_{M} \leq T^{\beta_{j+1}}} \frac{f(p_{1}...p_{M})}{\sqrt{p_{1}...p_{M}}} \nonumber \\
&& + T^{0.6}(\log\log T)^{2k} \nonumber \\
& \ll & T \exp\left(k^{2} \sum_{p \leq T^{\beta_{j}}} \frac{1}{p} \right) \cdot \frac{M!}{2^{M} (M/2)!} \left(\sum_{T^{\beta_{j}} < p \leq T^{\beta_{j+1}}} \frac{1}{p} \right)^{M/2} + T^{0.6}(\log\log T)^{2k} \nonumber \\
& \ll & T \exp\left(k^{2} \sum_{p \leq T^{\beta_{j}}} \frac{1}{p} \right) \cdot \left(\frac{1}{20\beta_{j+1}} \sum_{T^{\beta_{j}} < p \leq T^{\beta_{j+1}}} \frac{1}{p} \right)^{[1/(10\beta_{j+1})]} + T^{0.6}(\log\log T)^{2k} . \nonumber
\end{eqnarray}
Here the error term is $T^{0.6}(\log\log T)^{2k}$, rather than $T^{0.2}(\log\log T)^{2k}$ as in the proof of Lemma 1, because one picks up two additional factors from $\left( \Re G_{(j+1,l)}(t) \right)^M$ when handling the error term that comes from Proposition 2. Each factor has size at most $p_{1}...p_{M}$, where $T^{\beta_{j}} < p_{i} \leq T^{\beta_{j+1}}$ for all $i$, and that product is clearly at most $T^{M\beta_{j+1}} \leq T^{1/5}$. We also point out that our bound does not depend at all on $l$, because in the course of our arguments we upper bound the coefficients $(1/p^{1/2+1/(\beta_{l}\log T)}) \log(T^{\beta_{l}}/p)/\log(T^{\beta_{l}})$ of $G_{(j+1,l)}(t)$ by $1/\sqrt{p}$, which doesn't depend on $l$. (Note the step in the proof of Lemma 1 where $C(\tilde{p},\tilde{q})$ is replaced by $D(\tilde{p},\tilde{q})$.)

Putting together the foregoing calculations, remembering to include the sum over $l$ and the prefactor $(\beta_{j+1}^{3/2})^{[1/(10\beta_{j+1})]}$ from our initial manipulations, we conclude that
\begin{eqnarray}
&& \int_{t \in \mathcal{S}(j)} \exp\left( 2k \Re \sum_{p \leq T^{\beta_{j}}} \frac{1}{p^{1/2 + 1/(\beta_{j} \log T) + it}} \frac{\log(T^{\beta_{j}}/p)}{\log(T^{\beta_{j}})} \right) dt \nonumber \\
& \ll & (\mathcal{I}-j) T \exp\left(k^{2} \sum_{p \leq T^{\beta_{j}}} \frac{1}{p} \right) \cdot \left(\frac{\beta_{j+1}^{1/2}}{20} \sum_{T^{\beta_{j}} < p \leq T^{\beta_{j+1}}} \frac{1}{p} \right)^{[1/(10\beta_{j+1})]} . \nonumber
\end{eqnarray}
Now if $j=0$ then the left hand side is $\text{meas}(\mathcal{S}(0))$, whilst
$$ \mathcal{I} \leq \log\log\log T , \;\;\;\;\; \beta_{0} = 0, \;\;\;\;\; \beta_{1} = \frac{1}{(\log\log T)^{2}} , \;\;\;\;\; \sum_{p \leq T^{1/(\log\log T)^{2}}} \frac{1}{p} \leq \log\log T , $$
so we indeed obtain that $\text{meas}(\mathcal{S}(0)) \ll T e^{-(\log\log T)^{2}/10}$, as claimed in Lemma 2. If $1 \leq j \leq \mathcal{I}-1$ then we instead have
$$ \mathcal{I}-j \leq \frac{\log(1/\beta_{j})}{\log 20} , \;\;\;\;\; \sum_{T^{\beta_{j}} < p \leq T^{\beta_{j+1}}} \frac{1}{p} = \log \beta_{j+1} - \log \beta_{j} + o(1) = \log 20 + o(1) \leq 10 , $$
and therefore
\begin{eqnarray}
\int_{t \in \mathcal{S}(j)} \exp\left( 2k \Re \sum_{p \leq T^{\beta_{j}}} \frac{1}{p^{1/2 + 1/(\beta_{j} \log T) + it}} \frac{\log(T^{\beta_{j}}/p)}{\log(T^{\beta_{j}})} \right) dt & \ll & T \exp\left(k^{2} \sum_{p \leq T^{\beta_{j}}} \frac{1}{p} \right) e^{-\beta_{j+1}^{-1} \log(1/\beta_{j+1}) /21} \nonumber \\
& \ll & T \log^{k^{2}}T \cdot e^{-\beta_{j+1}^{-1} \log(1/\beta_{j+1}) /21} , \nonumber
\end{eqnarray}
as also claimed in Lemma 2.
\begin{flushright}
Q.E.D.
\end{flushright}

\subsection{Proof of Lemma 3}
To prove Lemma 3 we need to show, firstly, that
$$ \int_{t \in \mathcal{T}} \exp\left( 2k \Re \left(  \sum_{p \leq T^{\beta_{\mathcal{I}}}} \frac{1}{p^{1/2 + 1/(\beta_{\mathcal{I}} \log T) + it}} \frac{\log(T^{\beta_{\mathcal{I}}}/p)}{\log(T^{\beta_{\mathcal{I}}})} + \sum_{p \leq \log T} \frac{(1/2)}{p^{1 + 2it}} \right) \right) dt \ll_{k} T \log^{k^{2}}T . $$
We will sketch how to modify the proof of Lemma 1 to obtain this, and exactly similar modifications of the proof of Lemma 2 yield the other estimates claimed in Lemma 3.

For each $0 \leq m \leq (1/\log 2)\log\log T$, let us define $P_{m}(t) := \sum_{2^{m} < p \leq 2^{m+1}} \frac{(1/2)}{p^{1 + 2it}}$, and
$$ \mathcal{P}(m) := \{T \leq t \leq 2T : |\Re P_{m}(t)| > 2^{-m/10} , \; \text{but} \; |\Re P_{n}(t)| \leq 2^{-n/10} \; \forall m+1 \leq n \leq \frac{\log\log T}{\log 2} \}. $$
Clearly if $t$ belongs to none of the sets $\mathcal{P}(m)$, meaning that $|\Re P_{n}(t)| \leq 2^{-n/10}$ for all $n$, then we have $ \Re \sum_{p \leq \log T} \frac{(1/2)}{p^{1 + 2it}} = O(1)$, and so the part of $\int_{t \in \mathcal{T}}$ corresponding to such ``good'' $t$ can be bounded exactly as in Lemma 1. One can also easily show (e.g. by following the argument below but omitting the exponential factor) that $ \text{meas}(\mathcal{P}(m)) \ll T e^{-2^{3m/4}}$, and so if $2^{m} \geq (\log\log T)^{2}$, say, then the part of $\int_{t \in \mathcal{T}}$ corresponding to $t \in \mathcal{P}(m)$ is negligibly small by the Cauchy--Schwarz inequality. So we can restrict attention to those $t \in \mathcal{T} \cap \mathcal{P}(m)$ for $0 \leq m \leq (2/\log 2)\log\log\log T$.

Next, if $t \in \mathcal{P}(m)$ then we have
\begin{eqnarray}
\left| \Re \left(\sum_{p \leq 2^{m+1}} \frac{1}{p^{1/2 + 1/(\beta_{\mathcal{I}} \log T) + it}} \frac{\log(T^{\beta_{\mathcal{I}}}/p)}{\log(T^{\beta_{\mathcal{I}}})} + \sum_{p \leq \log T} \frac{(1/2)}{p^{1 + 2it}} \right) \right| & \leq & \sum_{p \leq 2^{m+1}} \frac{1}{\sqrt{p}} + \sum_{p \leq 2^{m+1}} \frac{(1/2)}{p} + O(1) \nonumber \\
& \ll & 2^{m/2} , \nonumber
\end{eqnarray}
and therefore we have
\begin{eqnarray}
&& \int_{t \in \mathcal{T} \cap \mathcal{P}(m)} \exp\left( 2k \Re \left(  \sum_{p \leq T^{\beta_{\mathcal{I}}}} \frac{1}{p^{1/2 + 1/(\beta_{\mathcal{I}} \log T) + it}} \frac{\log(T^{\beta_{\mathcal{I}}}/p)}{\log(T^{\beta_{\mathcal{I}}})} + \sum_{p \leq \log T} \frac{(1/2)}{p^{1 + 2it}} \right) \right) dt \nonumber \\
& \leq & e^{O(k2^{m/2})} \int_{t \in \mathcal{T} \cap \mathcal{P}(m)} \exp\left( 2k \Re \sum_{2^{m+1} < p \leq T^{\beta_{\mathcal{I}}}} \frac{1}{p^{1/2 + 1/(\beta_{\mathcal{I}} \log T) + it}} \frac{\log(T^{\beta_{\mathcal{I}}}/p)}{\log(T^{\beta_{\mathcal{I}}})} \right) dt \nonumber \\
& \leq & e^{O(k2^{m/2})} \int_{t \in \mathcal{T}} \left(2^{m/10} \Re P_{m}(t) \right)^{2[2^{3m/4}]} \exp\left( 2k \Re \sum_{2^{m+1} < p \leq T^{\beta_{\mathcal{I}}}} \frac{1}{p^{1/2 + 1/(\beta_{\mathcal{I}} \log T) + it}} \frac{\log(T^{\beta_{\mathcal{I}}}/p)}{\log(T^{\beta_{\mathcal{I}}})} \right) dt . \nonumber
\end{eqnarray}
Treating the exponential as in the proof of Lemma 1, noting $\Re P_{m}(t) = \sum_{2^{m} < p \leq 2^{m+1}} \frac{(1/2)\cos(2t\log p)}{p}$ and $\Re \sum_{2^{m+1} < p \leq T^{\beta_{\mathcal{I}}}} \frac{1}{p^{1/2 + 1/(\beta_{\mathcal{I}} \log T) + it}} \frac{\log(T^{\beta_{\mathcal{I}}}/p)}{\log(T^{\beta_{\mathcal{I}}})}$ are sums over disjoint sets of primes, we find that the above is
\begin{eqnarray}
& \ll & e^{O(k2^{m/2})} T \exp\left(k^{2} \sum_{2^{m+1} < p \leq T^{\beta_{\mathcal{I}}}} \frac{1}{p} \right) \cdot 2^{(m/10)2[2^{3m/4}]} \sum_{2^{m} < p_{1},...,p_{2[2^{3m/4}]} \leq 2^{m+1}} \frac{f(p_{1}...p_{2[2^{3m/4}]})}{p_{1}...p_{2[2^{3m/4}]}} + T^{0.2 + o(1)} \nonumber \\
& \ll & e^{O(k2^{m/2})} T \exp\left(k^{2} \sum_{2^{m+1} < p \leq T^{\beta_{\mathcal{I}}}} \frac{1}{p} \right) \cdot 2^{(m/10)2[2^{3m/4}]} \frac{(2[2^{3m/4}])!}{2^{2[2^{3m/4}]} ([2^{3m/4}])!} \left(\sum_{2^{m} < p \leq 2^{m+1}} \frac{1}{p^{2}} \right)^{[2^{3m/4}]} \nonumber \\
&& + T^{0.2+o(1)} .  \nonumber
\end{eqnarray}
Here we used the fact that if $2^{m} < p_{1},...,p_{2[2^{3m/4}]} \leq 2^{m+1}$ then $\prod p_{i} \leq 2^{(m+1)2[2^{3m/4}]} = T^{o(1)}$, to estimate the additional contribution to the error term (in Proposition 2) from $\left(\Re P_{m}(t) \right)^{2[2^{3m/4}]}$. (Recall that we have $2^{m} \leq \log T$.)

The above bound is
\begin{eqnarray}
& \ll & e^{O(k2^{m/2})} T \exp\left(k^{2} \sum_{2^{m+1} < p \leq T^{\beta_{\mathcal{I}}}} \frac{1}{p} \right) \cdot \left(2^{m/5} \cdot 2^{3m/4} \cdot \frac{1}{2^{m}} \right)^{[2^{3m/4}]} \nonumber \\
& \ll & e^{O(k2^{m/2}) - 2^{3m/4}} T \exp\left(k^{2} \sum_{2^{m+1} < p \leq T^{\beta_{\mathcal{I}}}} \frac{1}{p} \right) , \nonumber
\end{eqnarray}
which is $\ll e^{O(k2^{m/2}) - 2^{3m/4}} T \log^{k^{2}}T$. Then summing over $0 \leq m \leq (2/\log 2)\log\log\log T$ gives the desired bound $\ll_{k} T\log^{k^{2}}T$ for our original integral\footnote{The reader might wonder why we bothered to mention that we can restrict to values $m$ such that $2^{m} < (\log\log T)^{2}$, since this didn't appear to help us anywhere. When ``treating the exponential as in the proof of Lemma 1'' we need to have control over $|\Re \sum_{2^{m+1} < p \leq T^{\beta_{1}}} \frac{1}{p^{1/2 + 1/(\beta_{\mathcal{I}} \log T) + it}} \frac{\log(T^{\beta_{\mathcal{I}}}/p)}{\log(T^{\beta_{\mathcal{I}}})}|$, whereas the assumption that $t \in \mathcal{T}$ only implies that $|\Re \sum_{p \leq T^{\beta_{1}}} \frac{1}{p^{1/2 + 1/(\beta_{\mathcal{I}} \log T) + it}} \frac{\log(T^{\beta_{\mathcal{I}}}/p)}{\log(T^{\beta_{\mathcal{I}}})}| \leq \beta_{1}^{-3/4} = (\log\log T)^{3/2}$. But if $2^{m} < (\log\log T)^{2}$ then we trivially have that $\sum_{p \leq 2^{m+1}} \frac{1}{\sqrt{p}} = O(\log\log T)$.}.
\begin{flushright}
Q.E.D.
\end{flushright}

\section{Other moments of $L$-functions}
As promised in the introduction, in this section we indicate how the proof of Theorem 1 can be adapted to yield (presumably) sharp bounds for other moments of $L$-functions. Similarly as in Soundararajan's work~\cite{soundmoments}, all that one really needs is an analogue of Proposition 1 (upper bound by Dirichlet polynomials) and Proposition 2 (mean value/orthogonality result) for the $L$-functions under consideration. 

We will sketch the details in one important case, namely the moment of quadratic Dirichlet $L$-functions at the central point, where our method establishes the following:
\begin{thm2}
Assume the Generalised Riemann Hypothesis is true for all quadratic Dirichlet $L$-functions, and let $k \geq 0$ be fixed. Then for all large $X$ we have
$$ \sum_{\substack{X \leq |d| \leq 2X, \\ d \; \text{a fundamental discriminant}}} |L(1/2,\chi_{d})|^{k} \ll_{k} X\log^{k(k+1)/2}X , $$
where $\chi_{d} = \left(\frac{d}{\cdot}\right)$ denotes the primitive quadratic character corresponding to $d$.
\end{thm2}

Keating and Snaith~\cite{keatingsnaithL} conjectured an asymptotic formula for the left hand side in Theorem 2, building on work of Conrey and Farmer~\cite{conreyfarmer}, and the upper bound in Theorem 2 matches that conjecture. Unconditionally, Radziwi\l{\l} and Soundararajan~\cite{radsound} have proved a matching lower bound for all $k \geq 1$, and the method of Chandee and Li~\cite{chandeeli} would presumably yield a matching lower bound for all rational $0 \leq k < 1$ (although they do not explicitly discuss this example). See those papers for further references.

\vspace{12pt}
In this case one has the following results as analogues of Propositions 1 and 2.
\begin{prop3}
Let $X$ be large, and assume the Generalised Riemann Hypothesis is true for $L(s,\chi_{d})$ for all fundamental discriminants $d$ such that $X \leq |d| \leq 2X$. For each odd prime $p$, write $p* := (-1)^{(p-1)/2} p$. Then for any such $d$, and any $2 \leq x$, we have
\begin{eqnarray}
\log |L(1/2,\chi_{d})| & \leq & \sum_{p \leq x} \frac{\chi_{d}(p)}{p^{1/2+1/\log x}} \frac{\log(x/p)}{\log x} + (1/2)\log\log x + \frac{\log X}{\log x} + O(1) \nonumber \\
& = & \sum_{3 \leq p \leq x} \frac{\chi_{p*}(d)}{p^{1/2+1/\log x}} \frac{\log(x/p)}{\log x} + (1/2)\log\log x + \frac{\log X}{\log x} + O(1) . \nonumber
\end{eqnarray}
\end{prop3}

\begin{prop4}
Let $X$ be large and let $n=p_{1}^{\alpha_{1}}...p_{r}^{\alpha_{r}}$, where the $p_{i}$ are distinct odd primes and $\alpha_{i} \in \N$ for all $i$. Also write $p* = (-1)^{(p-1)/2} p$, as in Proposition 3. Then
$$ \sum_{X \leq |d| \leq 2X} \prod_{i=1}^{r} (\chi_{p_{i}*}(d))^{\alpha_{i}} \leq 2X \textbf{1}_{n \; \text{is a square}} + O(n) , $$
where $\textbf{1}$ denotes the indicator function.
\end{prop4}

The first inequality in Proposition 3 is essentially stated in $\S 4$ of Soundararajan's paper~\cite{soundmoments}, if one takes $\lambda = 1$ there (rather than $\lambda=\lambda_{0}$) and notes that the contribution from prime cubes and higher powers is $O(1)$. The second bound follows because the contribution from $p=2$ is clearly $O(1)$, and if $p$ is an odd prime and $d$ is a fundamental discriminant then $\chi_{d}(p) = \chi_{p*}(d)$, in view of quadratic reciprocity. See e.g. Theorem 9.14 in Montgomery and Vaughan's book~\cite{mv}, noting that $p*$ is always a fundamental discriminant, and that $\chi_{d}(p) = \chi_{d}(p*) \epsilon(d,p*)$ in the notation of that theorem.

One can choose any $x$ in Proposition 3, unlike in Proposition 1, because the $L$-function $L(s,\chi_{d})$ has no residue at 1 that risks producing any kind of main term. The fact that the contribution from prime squares in Proposition 3 is always of fixed size $\approx (1/2)\log\log x$, (because $\chi_{d}(p^{2}) = 1$ provided $p \nmid d$), rather than being given by a Dirichlet polynomial whose size can vary, reflects a genuine feature of this moment problem and will make the proof of Theorem 2 easier, since there will be no need to obtain an analogue of Lemma 3 to handle the prime squares contribution.

Proposition 4 follows trivially because $\prod_{i=1}^{r} \chi_{p_{i}*}^{\alpha_{i}}$ is a real character to modulus $n$, and is principal if and only if $n=p_{1}^{\alpha_{1}}...p_{r}^{\alpha_{r}}$ is a square.

\vspace{12pt}
Given Propositions 3 and 4, one can develop analogues of Lemmas 1 and 2 and use those to prove Theorem 2. More specifically, let us define
$$ \alpha_{0} := \frac{\log 2}{\log X}, \;\;\;\;\; \alpha_{i} := \frac{20^{i-1}}{(\log\log X)^{2}} \;\;\; \forall i \geq 1, $$
$$ \mathcal{J} = \mathcal{J}_{k,X} := 1 + \max\{i : \alpha_{i} \leq e^{-1000(1+k)}\} , $$
$$ \mathcal{Q} := \{X \leq |d| \leq 2X : \left| \sum_{X^{\alpha_{i-1}} < p \leq X^{\alpha_{i}}} \frac{\chi_{p*}(d)}{p^{1/2 + 1/(\alpha_{\mathcal{J}} \log X)}} \frac{\log(X^{\alpha_{\mathcal{J}}}/p)}{\log(X^{\alpha_{\mathcal{J}}})} \right| \leq \alpha_{i}^{-3/4} \; \forall 1 \leq i \leq \mathcal{J} \} . $$
Then one can show, analogously to Lemma 1, that
$$ \sum_{d \in \mathcal{Q}} \exp\left(k \sum_{3 \leq p \leq X^{\alpha_{J}}} \frac{\chi_{p*}(d)}{p^{1/2 + 1/(\alpha_{\mathcal{J}} \log X)}} \frac{\log(X^{\alpha_{\mathcal{J}}}/p)}{\log(X^{\alpha_{\mathcal{J}}})} \right) \ll X \log^{k^{2}/2}X , $$
and so
$$ \sum_{\substack{d \in \mathcal{Q}, \\ d \; \text{a fundamental discriminant}}} |L(1/2,\chi_{d})|^{k} \ll_{k} X \log^{k(k+1)/2}X , $$
in view of Proposition 3 with $x = X^{\alpha_{\mathcal{J}}}$. (Note how the factor $(1/2)\log\log x$ in Proposition 3 produces the additional factor $\log^{k/2}X$ on the right hand side.) One can also define obvious analogues of the sets $\mathcal{S}(j)$ from the proof of Theorem 1, and show as in Lemma 2 that the contribution to the moment from $d$ in those sets is small. This completes the proof of Theorem 2.

\vspace{12pt}
In $\S 4$ of Soundararajan's paper~\cite{soundmoments}, he also discusses the moment corresponding to all primitive characters to a given (prime) modulus, and the moment corresponding to quadratic twists of an elliptic curve $L$-function. The author expects that the proofs of Theorems 1 and 2 should adapt to those cases also, but has not checked all details.

\section{Remarks on the implicit constant}
The reader may check through the calculations in $\S 3$, and verify that the implicit constant we obtain in Theorem 1 is of the form $e^{e^{O(k)}}$ for large $k$ (and for $T$ large enough depending on $k$). Indeed, we have that $\int_{T}^{2T} |\zeta(1/2+it)|^{2k} dt$ is
\begin{eqnarray}
& \leq & \sum_{j=0}^{\mathcal{I}-1} \int_{t \in \mathcal{S}(j)} |\zeta(1/2+it)|^{2k} dt + \int_{t \in \mathcal{T}} |\zeta(1/2+it)|^{2k} dt \nonumber \\
& \leq & O_{k}\left(\sqrt{T e^{-(\log\log T)^{2}/10} \cdot T \log^{(2k)^{2}+1}T} \right) + \sum_{j=1}^{\mathcal{I}-1} e^{O(k^{3})} e^{-0.01k/\beta_{j}} T \log^{k^{2}}T + \int_{t \in \mathcal{T}} |\zeta(1/2+it)|^{2k} dt \nonumber \\
& \leq & O_{k}\left(\sqrt{T e^{-(\log\log T)^{2}/10} \cdot T \log^{(2k)^{2}+1}T} \right) + \sum_{j=1}^{\mathcal{I}-1} e^{O(k^{3})} e^{-0.01k/\beta_{j}} T \log^{k^{2}}T + e^{O(k^{3})} e^{2k/\beta_{\mathcal{I}}} T \log^{k^{2}}T , \nonumber
\end{eqnarray}
on inserting those calculations, and noting that the implicit constant in Lemma 3 is of the form $e^{O(k^{3})}$. (It is a sum over $m$ of terms of the form $e^{O(k2^{m/2}) - 2^{3m/4}}$.) The dominant term in the above is clearly the third one, as expected since it corresponds to most of the integral, and since $e^{-1000k} \leq \beta_{\mathcal{I}} \leq 20e^{-1000k}$ when $T$ is large we indeed find that
$$ \int_{T}^{2T} |\zeta(1/2+it)|^{2k} dt \ll e^{e^{O(k)}} T \log^{k^{2}} T , $$
with an absolute implied constant, when $k$ and $T$ are large.

This situation is a little disappointing, because the celebrated random matrix theory approach of Keating and Snaith~\cite{keatingsnaith}, now extensively developed, suggests that for any fixed $k \geq 0$ one should have
$$ \int_{0}^{T} |\zeta(1/2+it)|^{2k} dt \sim a(k) f(k) T \log^{k^{2}}T \;\;\;\;\; \text{as } T \rightarrow \infty , $$
where the ``arithmetical factor'' $a(k)$ is defined by
$$ a(k) = \prod_{p} \left(\left(1-\frac{1}{p}\right)^{k^{2}} \cdot \sum_{m=0}^{\infty} \left(\frac{\Gamma(m+k)}{m! \Gamma(k)}\right)^{2} \frac{1}{p^{m}} \right) , $$
and the other factor $f(k)$ is defined by
$$ f(k) = \lim_{N \rightarrow \infty} N^{-k^{2}} \prod_{j=1}^{N} \frac{\Gamma(j) \Gamma(j+2k)}{\Gamma(j+k)^{2}} . $$
In particular, one can show that $a(k) = e^{-k^{2}\log\log k + O(k^{2})}$ as $k \rightarrow \infty$, (the main contribution coming from primes $p \leq k^{2}$), and that $f(k) = e^{-k^{2}\log k + O(k^{2})}$, (the main contribution coming from $k \leq j \leq N$), so that $a(k)f(k) \leq e^{-k^{2}\log k}$ for large fixed $k$. See e.g. the introduction to Conrey and Gonek's paper~\cite{conreygonek}, and the papers of Conrey and Farmer~\cite{conreyfarmer} and of Diaconu, Goldfeld and Hoffstein~\cite{dgh}, for more discussion about the size of the constants in moment conjectures, from the random matrix and other points of view.

The basic reason that we obtain a poor constant in Theorem 1 is because we choose $\beta_{\mathcal{I}} \asymp e^{-1000k}$ rather small, so that the Dirichlet polynomial we work with in Lemma 1 is somewhat short. This is forced upon us by the step (\ref{sumlemma2}) in the proof of Theorem 1 where we sum the contribution coming from Lemma 2 over $j$. As discussed in the introduction, we look upon Lemmas 1, 2 and 3 as statements that certain integrals involving Dirichlet polynomials behave in the same way as random (e.g. Gaussian) models for those polynomials. Thus we might ask what conjectural result we would obtain by simply assuming that our Dirichlet polynomials behave in the same way as their Gaussian models, and then choosing the length of the polynomial in Proposition 1 optimally (and, in particular, choosing a longer polynomial than we can work with rigorously). It turns out that, with such a ``random Euler product'' model (that completely neglects the zeros of the zeta function, as in Proposition 1), one recovers roughly the same implicit constant, as an upper bound, as is implied by random matrix theory.

Indeed, for any $\log^{2}T \leq x \leq T^{2}$ we have that $\int_{T}^{2T} |\zeta(1/2+it)|^{2k} dt$ is
\begin{eqnarray}
= \int_{T}^{2T} e^{2k\log|\zeta(1/2+it)|} dt & \leq & e^{2k(\frac{\log T}{\log x} + O(1))} \int_{T}^{2T} e^{2k \Re\left(\sum_{p \leq x} \frac{1}{p^{1/2+1/\log x + it}} \frac{\log(x/p)}{\log x} + \sum_{p \leq \log T} \frac{(1/2)}{p^{1 + 2it}} \right)} dt \nonumber \\
& \approx & e^{2k(\frac{\log T}{\log x} + O(1))} T \E e^{2k N(0,(1/2)\log\log x + O(1))} , \nonumber
\end{eqnarray}
where the inequality follows from Proposition 1, and the final approximation is our heuristic that $\Re\left(\sum_{p \leq x} \frac{1}{p^{1/2+1/\log x + it}} \frac{\log(x/p)}{\log x} + \sum_{p \leq \log T} \frac{(1/2)}{p^{1 + 2it}} \right)$ should typically, as $T \leq t \leq 2T$ varies, behave in the same way as a Gaussian random variable with mean zero and variance $(1/2)\sum_{p \leq x} 1/p^{1+2/\log x} (\log(x/p)/\log x)^{2} + O(1) = (1/2)\log\log x + O(1)$. Then we have
\begin{eqnarray}
\E e^{2k N(0,(1/2)\log\log x + O(1))} & = & \frac{1}{\sqrt{\pi (\log\log x + O(1))}} \int_{-\infty}^{\infty} e^{2kz} e^{-z^{2}/(\log\log x + O(1))} dz \nonumber \\
& = & \frac{e^{k^{2}(\log\log x+O(1))}}{\sqrt{\pi (\log\log x + O(1))}} \int_{-\infty}^{\infty} e^{-(z- k(\log\log x + O(1)))^{2}/(\log\log x + O(1))} dz \nonumber \\
& = & e^{k^{2}(\log\log x+O(1))} , \nonumber
\end{eqnarray}
and so (conjecturally, for large $k$ and $T \gg_{k} 1$) we find $\int_{T}^{2T} |\zeta(1/2+it)|^{2k} dt$ is
$$ \lesssim \inf_{\log^{2}T \leq x \leq T^{2}} e^{2k(\frac{\log T}{\log x} + O(1))} T e^{k^{2}(\log\log x+O(1))} = T e^{O(k^{2})} \inf_{\log^{2}T \leq x \leq T^{2}} e^{2k\frac{\log T}{\log x} + k^{2}\log\log x} . $$
The infimum is attained when $\log x = (2\log T)/k$, yielding the bound $\int_{T}^{2T} |\zeta(1/2+it)|^{2k} dt \lesssim e^{-k^{2}\log k + O(k^{2})} T \log^{k^{2}}T$.

The reader might object that the above argument produces no arithmetical factor $a(k)$, but that is because we were only seeking an upper bound and didn't take any special care of the small prime contribution, which will not quite be Gaussian-like. Indeed, in the rigorous argument in $\S 4$ one already sees that, if one cares about the implicit constant, our treatment of the small primes gives quite a lot away. (Note the step after line (\ref{countingtuples}) where we upper bound $1/\prod_{j=1}^{r} \alpha_{j}!$ by 1, which is very generous if many of the $\alpha_{j}$ are large, as is the case for the small prime contribution.)

Thus a random Euler product model seems to produce very similar conclusions to random matrix theory models, but unfortunately we cannot make this rigorous, at the level of the implicit constant, without handling Dirichlet polynomials of length $T^{2/k}$ in Lemma 2, which the author is unable to do.

\vspace{12pt}
\noindent {\em Acknowledgements.} The author would like to thank Sandro Bettin, Andrew Granville, Maksym Radziwi\l{\l} and Kannan Soundararajan for some helpful remarks and discussions.

\end{document}